\documentclass[twocolumn]{autart}

\usepackage{cite}
\usepackage{xspace}
\usepackage{amsmath,amssymb,amsfonts}
\usepackage{algorithm}
\usepackage{algpseudocode}

\usepackage{graphicx}
\usepackage{textcomp}
\usepackage{xcolor}
\usepackage{url}
\usepackage{subfig}

\makeatletter
\DeclareRobustCommand{\qed}{%
	\ifmmode 
	\else \leavevmode\unskip\penalty9999 \hbox{}\nobreak\hfill
	\fi
	\quad\hbox{\qedsymbol}}
\newcommand{\openbox}{\leavevmode
	\hbox to.77778em{%
		\hfil\vrule
		\vbox to.675em{\hrule width.6em\vfil\hrule}%
		\vrule\hfil}}
\newcommand{\qedsymbol}{\openbox}
\newenvironment{proof}[1][\proofname]{\par
	\normalfont
	\topsep6\p@\@plus6\p@ \trivlist
	\item[\hskip\labelsep\itshape
	#1.]\ignorespaces
}{%
	\qed\endtrivlist
}
\newcommand{\proofname}{Proof}
\makeatother

\def\BibTeX{{\rm B\kern-.05em{\sc i\kern-.025em b}\kern-.08em
    T\kern-.1667em\lower.7ex\hbox{E}\kern-.125emX}}

\DeclareMathOperator{\ones}{\mathbf{1}}

\makeatletter
\newcommand\RedeclareMathOperator{%
	\@ifstar{\def\rmo@s{m}\rmo@redeclare}{\def\rmo@s{o}\rmo@redeclare}%
}
\newcommand\rmo@redeclare[2]{%
	\begingroup \escapechar\m@ne\xdef\@gtempa{{\string#1}}\endgroup
	\expandafter\@ifundefined\@gtempa
	{\@latex@error{\noexpand#1undefined}\@ehc}%
	\relax
	\expandafter\rmo@declmathop\rmo@s{#1}{#2}}
\newcommand\rmo@declmathop[3]{%
	\DeclareRobustCommand{#2}{\qopname\newmcodes@#1{#3}}%
}
\@onlypreamble\RedeclareMathOperator
\makeatother
\RedeclareMathOperator{\P}{\mathcal{P}}
\DeclareMathOperator{\N}{\mathcal{N}}
\DeclareMathOperator{\X}{\mathcal{X}}
\DeclareMathOperator{\Y}{\mathcal{Y}}
\DeclareMathOperator{\Pp}{\mathbb{P}}

\DeclareMathOperator{\B}{\mathcal{B}}

\DeclareMathOperator{\zeros}{\mathbf{0}}
\DeclareMathOperator{\SSS}{\mathcal{S}}
\DeclareMathOperator{\U}{\mathcal{U}}
\DeclareMathOperator{\W}{\mathcal{W}}

\newcommand{\real}{\ensuremath{\mathbb{R}}}
\newcommand{\complex}{\ensuremath{\mathbb{C}}}
\newcommand{\subscr}[2]{#1_{\textup{#2}}}
\newcommand{\supscr}[2]{#1^{\textup{#2}}}
\newcommand{\diag}[1]{\operatorname{diag}(#1)}
\newcommand{\until}[1]{\{1,\dots, #1\}}

\newcommand{\argmin}{\ensuremath{\operatorname{argmin}}}
\newcommand{\argmax}{\ensuremath{\operatorname{argmax}}}
\newcommand{\spn}{\operatorname{span}}
\newcommand{\image}{\operatorname{image}}
\newcommand{\setdef}[2]{\{#1 \; | \; #2\}}
\newcommand*{\QED}{\hfill\ensuremath{\square}}

\newcommand{\binpromc}{\textsc{Newton-like Neural Network}\xspace} 
\newcommand{\binpromd}{\textsc{Newton-like Neural Network}\xspace} 

\newcommand{\binpac}{\textsc{NNN-c}\xspace} 
\newcommand{\binpad}{\textsc{NNN-d}\xspace} 

\newcommand{\longthmtitle}[1]{\mbox{}{\bf \textit{(#1).}}}

\theoremstyle{definition}
\newtheorem{assump}{Assumption}

\begin{document}

	\begin{frontmatter}
		
		\title{Distributed Resource Allocation with Binary Decisions via
			Newton-like Neural Network Dynamics}
		
		\author{Tor Anderson\thanksref{footnoteinfo}}
		\author{\quad Sonia Mart{\'\i}nez\thanksref{footnoteinfo}}
		\thanks[footnoteinfo]{Tor Anderson
			and Sonia Mart{\'\i}nez are with the Department of Mechanical and
			Aerospace Engineering, University of California, San Diego, CA,
			USA. Email: {\small {\tt \{tka001, soniamd\}@eng.ucsd.edu}}. This
			research was supported by the Advanced Research Projects Agency -
			Energy under the NODES program, Cooperative Agreement
			DE-AR0000695.}


		\begin{abstract} This paper aims to solve a distributed resource allocation problem with binary local constraints.
			The problem is formulated as a binary program with a cost
			function defined by the summation of
			agent costs plus a global mismatch/penalty
			term.
			We propose a
			modification of the Hopfield Neural Network (HNN) dynamics in order to
			solve this problem while incorporating a novel Newton-like weighting
			factor. This addition lends itself to fast avoidance of saddle
			points, which the gradient-like HNN is susceptible to. Turning to a
			multi-agent setting, we reformulate the problem and develop a
			distributed implementation of the Newton-like dynamics. We show that
			if a local solution to the distributed reformulation is obtained, it
			is also a local solution to the centralized problem. A main
			contribution of this work is to show that the probability of
			converging to a saddle point of an appropriately defined energy
			function in both the centralized and distributed settings is zero
			under light
			assumptions.
			Finally, we enlarge our algorithm with an annealing technique which
			gradually learns a feasible binary solution. Simulation results
			demonstrate that the proposed methods are competitive with
			centralized greedy and SDP relaxation approaches in terms of
			solution quality, while the main advantage of our approach is a
			significant improvement in runtime over the SDP relaxation method
			and the distributed quality of implementation.
		\end{abstract}
		
	\end{frontmatter}

\section{Introduction}

There has been an explosion of literature
surrounding the design of distributed algorithms for convex
optimization problems and how these pertain to the operation of future
power grids. 
A common assumption of these algorithms is the property of
\emph{convexity}, which lends itself to provably optimal solutions
which are \emph{scalable} and \emph{fast}. However, some settings give rise to nonconvex decision sets. For example, in an optimal power dispatch setting, devices available
for providing load-side frequency regulation such as HVAC systems, household appliances, and manufacturing
systems are often limited to
discrete on/off operational modes. It is even preferable to charge populations of electric
vehicles in a discrete on/off manner due to nonlinear battery chemistries. The available tools in optimization for these
nonconvex settings are less mature, and when considering a distributed
setting in which devices act as agents that collectively compute a
solution over a sparse communication graph, the available tools are
significantly less developed.  With this in mind, we
are motivated to develop a \emph{scalable}, \emph{fast} approach for
these binary settings which is amenable to a \emph{distributed}
implementation.

Quadratic programs with nonconvex binary
constraints are known to be NP-hard in general,
see~\cite{PC-AS:95,DL-XS-SG-JG-CL:10}. In this paper, we consider a problem which is quite applicable to the economic dispatch problem in power networks, see~\cite{XH-XF-JY:19,XH-JY-TH-CL:19,BH-LL-HZ-YL-QS:19} for three recent examples in microgrid environments. However, none of these examples address devices with binary constraint sets. The binary problem is, however, desirable to approach in a distributed context~\cite{ZY-AB-HZ-NZ-QX-CK:16,PY-YH-LF:16}. Greedy
algorithms~\cite{TC-CL-RR-CS:09} have been proposed for binary
programs, such as the well-known Traveling Salesman Problem (TSP), but
it is well documented that these methods can greatly suffer in
performance~\cite{GG-AY-AZ:02} except in cases where the cost function
is submodular~\cite{GLN-LAW-MLF:78,MS-SB-HV:10}. A more modern
approach to solving optimization problems with a binary feasibility
set is to cast them as a semidefinite program (SDP) with a nonlinear
rank constraint, see~\cite{SP-FR-HW:95,LV-SB:96,SB-LV:97} for some
classical references or~\cite{ZQL-WKM-AMCS-YY-SZ:10,PW-CS-AH-PT:17}
for more recent work on the topic.
By relaxing the rank constraint, a convex problem is obtained whose
solution can be shown to be equal to the optimal dual value of the
original problem, see e.g.~\cite{PP-SL:03}. However, it is necessary
in these approaches to either impose a single centralized coordinator
to compute the solution and broadcast it to the actuators or agents,
or schedule computations, which suffers from scalability issues,
privacy concerns, and does not enjoy the simpler and more robust
implementation of a distributed architecture in a large network.

Neuro-dynamic programming is a different paradigm for addressing
nonconvex problems with computational tractability,
see~\cite{DPB-JNT:96} for a broad reference. A neural-network based
method for binary programs was first developed by Hopfield
in~\cite{JH-DT:85}, which was originally proposed in order to address
TSPs. We refer to this method from here on as a Hopfield Neural
Network (HNN). This method provided a completely different avenue for
approaching binary optimizations, and followup works are found
in~\cite{KS:96,JM:96,BKP-BKP:92,SB-ZA:02}. These works formalize and
expand the framework in which the HNN method is applicable. However,
these algorithms essentially implement a gradient-descent on an
applicable nonconvex energy function, which is susceptible to being
slowed down by convergence to saddle-points. There are avenues for
Newton-like algorithms in nonconvex environments to address this
issue, which incorporate some treatment of the negative Hessian
eigenvalues in order to maintain a monotonic descent of the cost
function, see e.g.~\cite{PG-WM-MW:81,YD-RP-CG-KC-SG-YB:14}. A recently
developed method employs a Positive-definite Truncated inverse
(PT-inverse) operation on the Hessian of a nonconvex energy or cost
function in order to define a nonconvex Newton-descent
direction~\cite{SP-AM-AR:19}, although the technique does not
presently address binary settings. Perhaps more importantly, all
variants of existing HNN methods and the aforementioned works for
nonconvex Newton-like algorithms are framed for centralized
environments in which each agent knows global information about the
state of all other agents, which is not
scalable.

The contributions of this paper are threefold. We start by considering
a binary programming problem formulated as a summation of local costs
plus a squared global term. By leveraging a specific choice for the
cost functions, we adapt the setting to an HNN framework. Then, we
propose a novel modification of the dynamics with a PT-inverse of the
Hessian of an appropriate energy function to define centralized
$\binpromc$ ($\binpac$). We prove a rigorous convergence result 
to a local minimizer, thus excluding saddle-points, with probability
one, given some mild assumptions on the algorithm parameters and
initial condition.
Thirdly, we reformulate the problem so that it is solvable via a
distributed algorithm by means of an auxiliary variable. We show that
local solutions of the distributed reformulation are equivalent to
local solutions of the centralized one, and we define a corresponding
energy function and distributed algorithm for which we show
convergence to a local minimizer with probability
one. 
Simulations validate that our method is
superior to SDP relaxation approaches in terms of runtime and
scalability and outperforms greedy methods in terms of
scalability.

\section{Preliminaries}

This
section establishes notation\footnote{The set of real numbers, real positive numbers, real
  $n$-dimensional vectors, and real $n$-by-$m$ matrices are written as
  $\real,\real_+, \real^n$, and $\real^{n\times m}$, respectively. We denote by
  $x_i$ the $\supscr{i}{th}$ element of $x\in\real^n$ and $A_{ij}$ the
  element in the $\supscr{i}{th}$ row and $\supscr{j}{th}$ column of
  $A\in\real^{n\times m}$. 
  For a square matrix $A$, we denote by $A^{\dagger}$ the
  Moore-Penrose pseudoinverse of $A$. We use the shorthand $\ones_n =
  (1,\dots,1)^\top \in\real^n$ and $\zeros_n = (0,\dots,0)^\top
  \in\real^n$. Cartesian products of sets are denoted by a
  superscript, for example, $\{0,1\}^n =
  \{0,1\}\times\dots\times\{0,1\}$. The gradient of a function
  $f:\real^{n}\rightarrow\real$ with respect to $x\in \real^n$ at $x$
  is denoted by $\nabla_x f(x)\in \real^n$, and the Hessian matrix of
  $f$ at $x$ is written as $\nabla_{xx} f(x)\in\real^{n\times n}$. We denote elementwise operations on vectors $x,y\in\real^n$ as $(x_i y_i)_i = (x_1 y_1, \dots, x_n y_n)^\top$, $(x_i)_i^2 = (x_1^2, \dots, x_n^2)^\top$,
  $(c/x_i)_i = (c/x_1, \dots, c/x_n)^\top$, $\log (x_i)_i = (\log(x_1), \dots, \log(x_n))^\top$, and $(e^{x_i})_i = (e^{x_1}, \dots, e^{x_n})^\top$. The notation $\diag{x}$ indicates the diagonal matrix with entries given by elements of $x$, and
  $\B(x,\eta )$ denotes the closed ball of radius $\eta$ centered at
  $x$.} 
and background concepts to be used throughout
the paper.

We refer the reader to~\cite{FB-JC-SM:09} as a Graph Theory supplement. One can define a
Laplacian matrix $L$ associated with a graph $\mathcal{G}$ as follows:
\begin{equation*}
L_{ij} = \begin{cases}
-1, & j\in\N_i, \\
-\sum_{k\neq i} L_{ik}, & i = j, \\
0, & \text{otherwise,}
\end{cases}
\end{equation*}
where $\N_i$ is the set of neighbors of node $i$. An immediate property is that $0$ is an eigenvalue of $L$
associated with the eigenvector $\ones_n$. It is simple iff $\mathcal{G}$ is
connected.

Next, we introduce the Positive-definite Truncated inverse
(PT-inverse) and its relevance to nonconvex Newton methods.

\begin{defn}[\hspace{1sp}\cite{SP-AM-AR:19}]\longthmtitle{PT-inverse}\label{def:ptinv}
  Let $A\in\real^{n\times n}$ be a symmetric matrix with an
  orthonormal basis of eigenvectors $Q\in\real^{n\times n}$ and
  diagonal matrix of eigenvalues $\Lambda\in\real^{n\times n}$. Consider a constant $m>0$ and define
  $\vert\Lambda\vert_m\in\real^{n\times n}$ by:
	\begin{equation*}
          (\vert \Lambda \vert_m)_{ii} = \begin{cases} \vert\Lambda_{ii}\vert, & \vert\Lambda_{ii}\vert \geq m, \\
            m, & \text{otherwise.}
	\end{cases}
	\end{equation*}
	The PT-inverse of $A$ with parameter $m$ is defined by $(\vert
        A\vert_m)^{-1} = Q^\top (\vert\Lambda\vert_m)^{-1} Q \succ 0$.
\end{defn}

The PT-inverse operation flips the sign on the negative eigenvalues of
$A$ and truncates near-zero eigenvalues to a (small) positive value
$m$ before conducting the inverse. Effectively, this generates a
positive definite matrix bounded away from zero to be inverted,
circumventing near-singular cases. In terms of
  computational complexity, it is on the order of standard
  eigendecomposition (or more generally, singular value
  decomposition), which is roughly $O(n^3)$~\cite{VP-ZC:99}. However,
  we note in Section~\ref{sec:dist-hop} that the matrix to be
  PT-inverted is diagonal, which is $O(n)$.

The PT-inverse is useful for nonconvex Newton
approaches~\cite{SP-AM-AR:19} in the following sense: first, recall
that the Newton descent direction of $f$ at $x$ is computed as
$-\left(\nabla_{xx} f(x)\right)^{-1}\nabla_x f(x)$. For $f$ strictly
convex, it holds that $\nabla_{xx} f(x) \succ 0$ and the Newton
direction is well defined and decreases the cost. For (non-strictly)
convex or nonconvex cases, $\nabla_{xx}f(x)$ will be singular,
indefinite, or negative definite. A PT-inverse operation remedies
these cases and preserves the descent quality of the
method. Additionally, saddle points are a primary concern for
first-order methods in nonconvex settings~\cite{YD-RP-CG-KC-SG-YB:14},
and the Newton flavor endowed by the PT-inverse effectively performs a
change of coordinates on saddles with ``slow" unstable manifolds
compared to the stable manifolds. We discuss this further in
Section~\ref{sec:all-hop}.

\section{Problem Statement and Dual Problem}

Here, we formally state the nonconvex optimization problem
we wish to solve and formulate its dual for the sake of deriving a
lower bound to the optimal cost.

We aim to find an adequate solution to a resource
  allocation problem where the optimization variables take the form
of binary decisions over a population of $n$ agents. We
  note that the problem we consider is  applicable to generator
  dispatch and active device response in an economic dispatch power
  systems setting~\cite{CAISO-BPM:18}, but the remainder of the paper
  will frame it primarily as resource allocation. Let each
agent $i\in\until{n}$ be endowed with a
  decision variable $x_i$ and a cost $c_i\in\real$, a value which
  indicates the incremental cost of operating in the $x_i=1$ state
  versus the $x_i=0$ state. We do not impose a sign restriction on
  $c_i$, but this may be a common choice in the power systems setting
  where $x_i=1$ represents an ``on" device state and $x_i=0$
  represents ``off." Additionally, each agent is endowed with a
  parameter $p_i$ which represents some incremental consumption or
  generation quantity when operating in the $x_i=1$ state versus
  $x_i=0$ and also a passive cost $d_i$.

We are afforded some design choice in the cost function
models for $x_i\notin\{0,1\}$, and for each $i \in
  \until{n}$, so we design abstracted cost functions
$f_i:[0,1]\rightarrow \real$ that satisfy $f_i(0) = d_i$
and $f_i(1) = c_i + d_i, \forall i$. This
  design choice is intrinsic to a cost model for any separable binary
  decision optimization context. In particular, the value of
  $f_i(x_i)$ for any $x_i\notin\{0,1\}$ is only relevant to the
  algorithm design, but need not have a physical interpretation or
  pertain to the optimization model since these points are
  infeasible. With this in mind, we enlarge the cost
  model by adopting the following:
\begin{assump}\longthmtitle{Quadratic Cost Functions}\label{ass:costs}
	The local cost functions $f_i$ take the form
	\begin{equation*}
	f_i(x_i) = \dfrac{a_i}{2} (x_i - b_i)^2 - \dfrac{a_ib_i^2}{2} + d_i,
	\end{equation*} with $a_i, b_i, d_i \in \real$.
\end{assump}
Note that, for any value $c_i = f_i(1)-f_i(0)$, there exists a
family of coefficients $a_i, b_i$ such that $(a_i/2) (1 -
b_i)^2 - (a_i/2) b_i^2 = c_i$. Further, the constant terms ensure $f_i(0) = d_i$ and $f_i(1) = c_i + d_i$. The design of
$a_i, b_i$ will be discussed in Section~\ref{sec:all-hop}.

The problem we aim to solve can now be formulated as:
\begin{equation*} 
  \P 1: \
  \underset{x\in\{0,1\}^n}{\text{min}} \
  f(x) = \sum_i^n f_i(x_i) + \dfrac{\gamma}{2}\left(p^\top x - \subscr{P}{r}\right)^2.  
\end{equation*}
Here, $\subscr{P}{r}\in\real$ is a given reference value to be matched
by the total output $p^\top x$ of the devices, with
$p\in\real^n$ having entries $p_i$. This matching is enforced by means
of a penalty term with coefficient $\gamma >0$ in $\P 1$.
In the power systems setting, $\subscr{P}{r}$ can
  represent a real-power quantity to be approximately matched by the
  collective device-response. The coefficient $\gamma$ and the signal
  $\subscr{P}{r}$ are determined by an Independent System Operator
(ISO) and communicated to a Distributed Energy Resource Provider
(DERP) that solves $\P 1$ to obtain a real-time dispatch solution,
see~\cite{CAISO-BPM:18} for additional information.

  The primal $\P 1$ has an associated dual $\mathcal{D}1$ which takes
  the form of a semidefinite program (SDP) whose optimal value lower
  bounds the cost of $\P 1$. This SDP is
\begin{subequations} \label{eq:dual-prob}
	\begin{align}
	\mathcal{D}1: \ & \underset{\mu\in\real^n,\Delta\in\real}{\text{max}}
	& & \Delta,  \label{eq:dual-cost}\\
	& \text{subject to} & & \begin{bmatrix} \dfrac{1}{2}Q(\mu) &
	\xi(\mu) \\ \xi(\mu)^\top & \zeta - \Delta
	\end{bmatrix} \succeq 0. \label{eq:gen-fiedler-const}
	\end{align}
\end{subequations}
In $\mathcal{D}1$, $Q : \real^n \rightarrow \real^{n\times n}$ and
$\xi: \real^n \rightarrow \real^n$ are real-affine functions of $\mu$ and $\zeta$ is a constant. These definitions are $Q(\mu) = \left(\diag{a/2 + \mu} + \dfrac{\gamma}{2}pp^\top\right),
\xi(\mu) = ((a_i b_i)_i + \mu + \gamma \subscr{P}{r}p),$ and
$\zeta = \sum_{i=1}^n \dfrac{a_i b_i^2}{2} + \dfrac{\gamma}{2}\subscr{P}{r}^2$. See~\cite{SB-LV:04} for more detail on the derivation of $\mathcal{D}1$. 

\section{Centralized Newton-like Neural Network}\label{sec:all-hop}

In this section, we develop the Centralized\\$\binpromc$, or 
$\binpac$, which is well suited for solving $\P1$ in a centralized
 setting.

To draw analogy with the classic Hopfield Neural Network approach we
will briefly introduce an auxiliary variable $u_i$ whose relation to
$x_i$ is given by the logistic function $g$ for each $i$:
\begin{align}
\label{eq:act-fn}
x_i &= g(u_i) = \dfrac{1}{1+e^{-u_i/T}}, \quad &&u_i\in\real, \\
u_i &= g^{-1}(x_i) = -T \ \log\left(\dfrac{1}{x_i} - 1\right), \quad
&&x_i\in\left( 0, 1\right), \nonumber
\end{align}
with temperature parameter $T > 0$. 

Let $x\in (0,1)^n, u\in\real^n$ be vectors with entries given by $x_i,
u_i$. To establish our algorithm, it is appropriate to first define an
energy function related to $\P 1$. Consider
\begin{equation}\label{eq:energyc}
E(x) = f(x) + \dfrac{1}{\tau}\sum_i \int_{0}^{x_i} g^{-1} (\nu) d\nu,
\end{equation}
where $\tau > 0$ is a time-constant and for $z \in [0,1]$,
\vspace{-0.5cm}{\small
	\begin{equation*}
	\int_{0}^{z} g^{-1}(\nu)d\nu = \begin{cases}
	T\left(\log(1-z) - z \log (\frac{1}{z} - 1)\right), & z\in(0,1), \\
	0, & z\in\{0,1\}. 
	\end{cases}
      \end{equation*}}
    The classic HNN implements dynamics of the form $\dot{u} = -\nabla_x
    E(x)$, where the equivalent dynamics in $x$
    can be computed as $\dot{x} = - \nabla_x E(x)dx/du$. 
     These dynamics can  be thought of to
      model the interactions between neurons in a neural network or
      the interconnection of amplifiers in an electronic circuit,
      where in both cases the physical system tends toward low energy
      states, see~\cite{JH-DT:85,KS:96}. 
      In an optimization setting, low energy 
      states draw analogy to low cost solutions. We 
      now describe our modification to the classical HNN dynamics.

Recall that the domain of $x$ is $(0,1)^n$ and our elementwise notation for $\log$ and division. We have the expressions $\nabla_x E(x) = -Wx -v - (T/\tau) \log\left(1/x_i-1\right)_i$ and $dx/du = (x-(x_i^2)_i)/T$,
where $W = -\diag{a} - \gamma pp^\top\in \real^{n\times n}$ and $v = (a_i b_i)_i + \gamma P_{r} p\in \real^n$ are defined via $f$. From this point forward, we work mostly in terms of $x$ for the sake
of consistency.  Consider modifying the classic HNN dynamics with a
PT-inverse $(\vert H(x) \vert_m)^{-1} \succ 0$ as in~\cite{SP-AM-AR:19}, where $H(x) = \nabla_{xx} E(x)$. The $\binpac$ dynamics are then given~by:
\begin{equation}\label{eq:binpac}
\begin{aligned}
\dot{x} &= -(\vert H(x)\vert_m)^{-1} \diag{\frac{dx}{du}} \nabla_x{E(x)} \\
&=(\vert H(x)\vert_m)^{-1} \diag{\frac{(x_i-x_i^2)_i}{T}} \\
&\qquad \qquad \qquad \left(Wx + v +
\dfrac{T}{\tau}\log\left(1/x_i-1\right)_i\right).
\end{aligned}
\end{equation}

These dynamics lend
to the avoidance of saddle points of $E$. To see this, consider
the eigendecomposition $H(\tilde{x}) = Q^\top \Lambda Q$ at some
$\tilde{x}$ near a saddle point, i.e. $\nabla_x E(\tilde{x})\approx
0$. If many entries of $\Lambda$ are small in magnitude and remain
small in the proximity of $\tilde{x}$, then the gradient is
changing slowly along the ``slow" manifolds associated with the eigenspace of
the small eigenvalues. This is precisely what the PT-inverse is
designed to combat: the weighting of the
dynamics is increased along these manifolds by a factor that is
inversely proportional to the magnitude of the
eigenvalues. Additionally, negative eigenvalues of the Hessian
are flipped in sign, which causes attractive
manifolds around saddle points to become repellent.

It is desirable for $E$ to be concave on most of its domain so the
trajectories are pushed towards the feasible points of $\P 1$;
namely, the corners of the unit hypercube. To examine this, the Hessian of $E$ can be
computed as $H(x) = \frac{d^2 f}{dx^2} + \frac{1}{\tau}\diag{\frac{dg^{-1}(x)}{dx}} = -W + \dfrac{T}{\tau}\diag{\frac{1}{(x_i-x_i^2)_i}}.$
Notice that the second term is positive definite on $x\in(0,1)^n$ and
promotes the convexity of $E$, particularly for elements $x_i$ close to $0$ or $1$. For a fixed $T,\tau$, choosing $a_i <
-\gamma\Vert p\Vert^2 - 4T/\tau, \forall i$ guarantees $E(x) \prec 0$
at $x=(0.5) \ones_n$. 
Generally speaking, choosing $a_i$ to be negative and
large in magnitude lends itself to concavity of $E$ over a larger subset of its domain and to trajectories converging closer
to the set $\{0,1\}^n$. However, this comes at the expense of not exploring a rich subset
of the domain. At the end of this section, we develop a Deterministic
Annealing (DA) approach inspired by~\cite{KR:98} for the online
adjustment of $T,\tau$ to obtain an effective compromise between
exploration of the state space and convergence to a feasible point of
$\P1$.

We now characterize the equilibria of~\eqref{eq:binpac} for
$x\in[0,1]^n$. 
It would appear that $x$ with some components $x_i\in\{0,1\}$ are
candidate equilibria due to the $x_i-x_i^2$ factor vanishing. However,
the dynamics are not well defined here due to the $\log$ term. Additionally, note that $\lim_{x_i\rightarrow \delta} e_i^\top H(x) e_i = \infty, \ \delta \in \{0,1\}, \forall i,$
where $e_i$ is the $\supscr{i}{th}$ canonical basis vector. Due to the
$\frac{T}{\tau(x_i - x_i^2)}$ term dominating $W$ in the expression for $H$ when $x_i$ values
are close to $\{0,1\}$, it follows that an eigenvalue of $(\vert
H(x)\vert_m)^{-1}$ approaches zero as $x_i\rightarrow 0$ or $1$ with
corresponding eigenvector approaching $v_i = e_i$:
\vspace{-0.5cm}{\small{\begin{equation*}
  \lim_{x_i\rightarrow \delta} = v_i^\top (\vert H(x)\vert_m )^{-1}
  v_i = \frac{T}{\tau}(x_i - x_i^2) = 0, \quad 
  \delta \in \{0,1\}, \forall i.
\end{equation*}}}
Using this fact, and ignoring $T,\tau > 0$, we can compute the 
undetermined limits in the components of $\dot{x}$ as $x_i \rightarrow
\delta \in \{0,1\}$ by repeated applications of L'Hospital's rule:
\begin{equation}\label{eq:lhospital-eq}
  \lim_{x_i\rightarrow \delta} \log\left(\dfrac{1}{x_i}-1\right)(x_i - x_i^2)^2 
  = \begin{cases}
  0, & \delta = 0^+, \\
  0, & \delta = 1^-. \end{cases}
\end{equation}
Thus, components $x_i\in\{0,1\}$ constitute candidate equilibria. We
will, however, return to the first line of~\eqref{eq:lhospital-eq} in the proof of
Lemma~\ref{lem:fwd-inv-c} to show that they are unstable.  As for
components of $x$ in the interior of the hypercube, the expression
$\dot{x} =0$ can not be solved for in closed form. However, we provide
the following Lemma which shows that the set of equilibria is finite.
\begin{lem}\longthmtitle{Finite Equilibria}~\label{lem:finite-eq}
  Let $\X$ be the set of equilibria of~\eqref{eq:binpac}
  satisfying $\dot{x} = 0$ on $x\in[0,1]^n$. The set $\X$ is finite.
\end{lem}
The proof can be found in the Appendix, and all proofs for the remainder of the paper will be contained there. 

\begin{figure}[h]
  \centering
	\includegraphics[scale = 0.4]{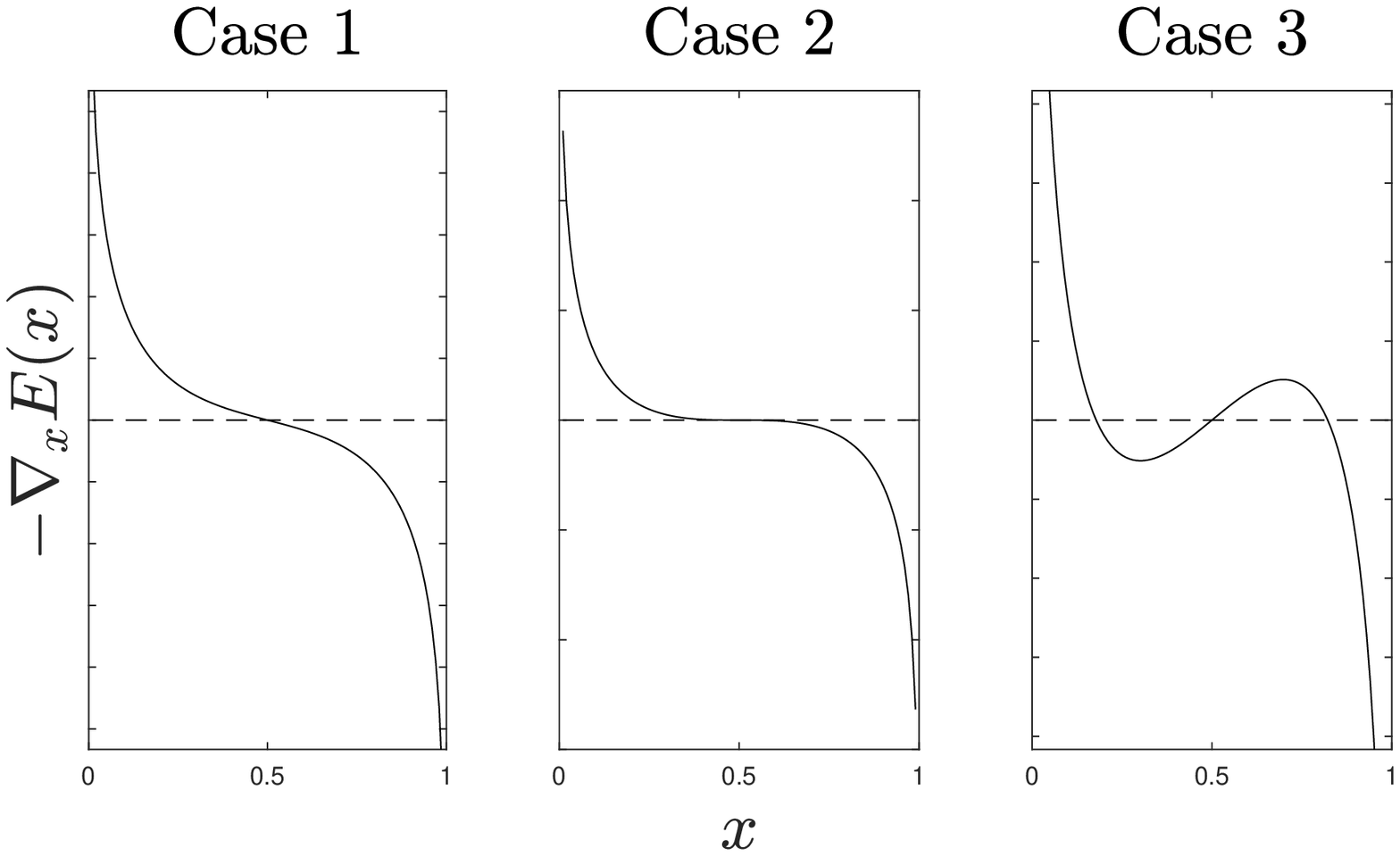}
	\includegraphics[scale = 0.4]{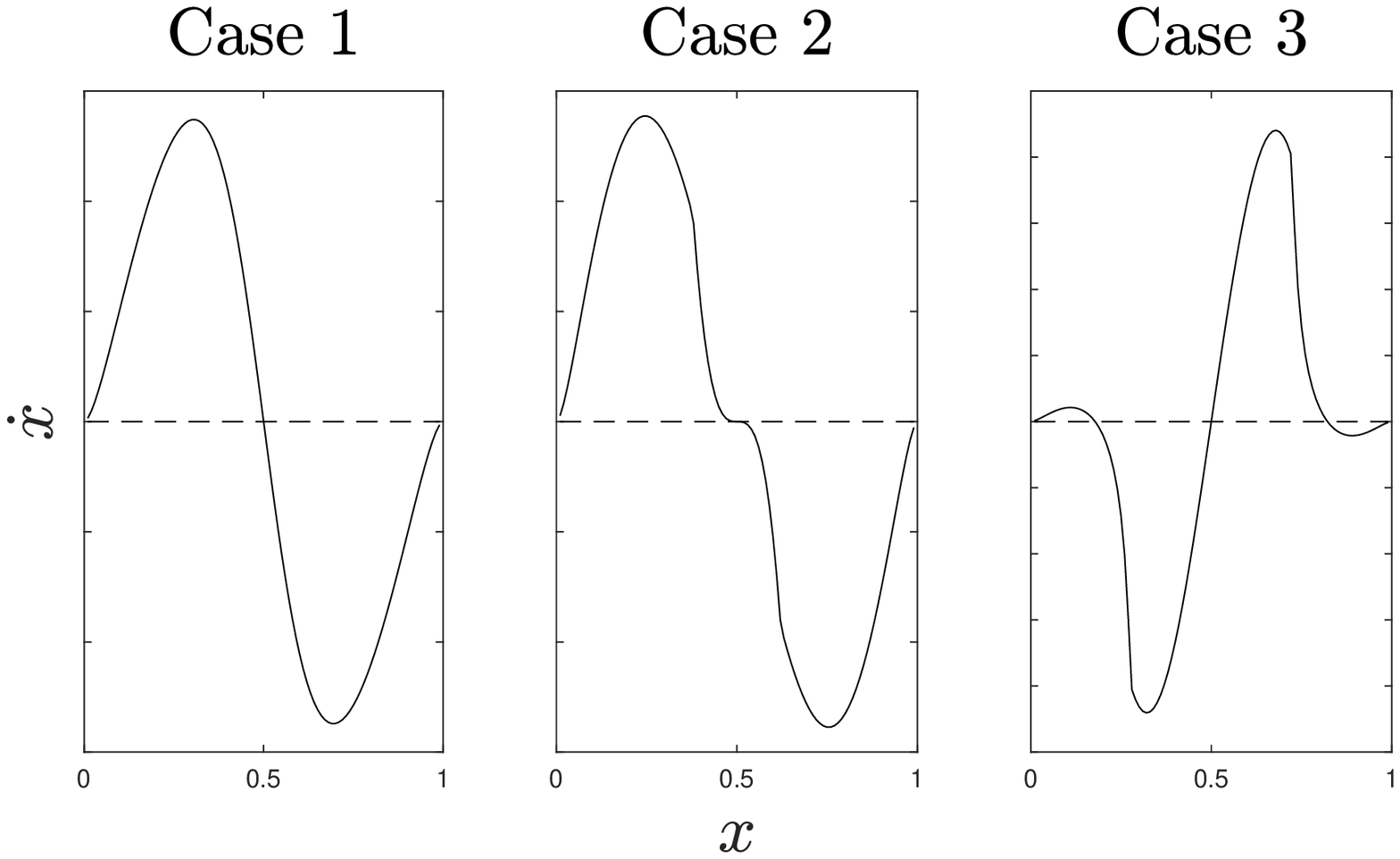}
	\caption{Illustration of $-\nabla_x E(x)$ (top) and $\dot{x}$
          (bottom) for three instances of $a$. Case 1: $a > -\gamma
          \Vert p \Vert ^2 - 4T/\tau$, Case 2: $a = -\gamma \Vert p
          \Vert ^2 - 4T/\tau$, Case 3: $a < -\gamma \Vert p \Vert ^2 -
          4T/\tau$.}
	\label{fig:xdot1d}
\end{figure}

To demonstrate the qualitative behavior of equilibria in a simple
case, consider a one-dimensional example with $a > - \gamma p^2 -
4T/\tau$ and recall that, for $x\in(0,1)$, the sign of $-\nabla_x
E(x)$ is the same as $\dot{x}$. In Figure~\ref{fig:xdot1d}, we observe
that $-\nabla_x E(x)$ monotonically decreases in $x$, and a globally
stable equilibrium exists in the interior $x\in(0,1)$ near
$x=0.5$. On the other hand, $a < - \gamma p^2 - 4T/\tau$ gives way
to $3$ isolated equilibria in the interior (one locally unstable near
$x=0.5$ and two locally stable near $x\in\{0,1\}$). This behavior
extends in some sense to the higher-dimensional case. Therefore, for a
scheme in which $T$ and $\tau$ are held fixed, we prescribe $a < -
\gamma \| p\|^2 - 4T/\tau$. We provide a Deterministic Annealing
(DA) approach inspired by~\cite{KR:98} for the online adjustment of $T,\tau$ in the following subsection which
compromises with this strict design of $a$.

Finally, we establish a Lemma about the domain of the trajectories
of~\eqref{eq:binpac}.
\begin{lem}\longthmtitle{Forward Invariance of the Open Hypercube}\label{lem:fwd-inv-c}
  The open hypercube $(0,1)^n$ is a forward-invariant set under the
  $\binpac$ dynamics~\eqref{eq:binpac}.
\end{lem}

Knowing that $\P 1$ is generally NP-hard, it is unlikely
that a non-brute-force algorithm exists that can converge to a global
minimizer. For this reason, we aim to establish
asymptotic stability to a local minimizer of $E$. We first establish some assumptions.

\begin{assump}\longthmtitle{Random Initial Condition}\label{ass:init}
  The initial condition $x(0)$ is chosen randomly according to a
  distribution $\Pp$ that is nonzero on sets that have nonzero volume
  in $[0,1]^n$.
\end{assump}
An appropriately unbiased initial condition for our algorithm is
$x(0)\approx(0.5)\ones_n$, 
which is adequately far from the local minima located near corners of
the unit cube. So, we suggest choosing a uniformly random $x(0) \in \B
((0.5)\ones_n ,\epsilon)$, where $0 < \epsilon \ll
1$.

\begin{assump}\longthmtitle{Choice of $T,\tau$}\label{ass:t-tau}
  The constants $T, \tau > 0$ are each chosen randomly according to a
  distribution $\bar{\Pp}$ that is nonzero on sets that have nonzero
  volume on $\real_+$. 
\end{assump}
Similarly to $x(0)$, we suggest choosing these constants uniformly
randomly in a ball around some nominal $T_0, \tau_0$,
i.e. $T\in\B(T_0,\epsilon), \tau\in\B(\tau_0,\epsilon), 0 < \epsilon
\ll 1$. The $T_0,\tau_0$ themselves are design parameters stemming from the neural network model, and we provide some intuition for selecting these in the simulation Section.

Now we state the main convergence result of $\binpac$ in Theorem~\ref{thm:cent-cvg}, which states that for a random choice of $T,\tau$, an
	initial condition chosen randomly from $(0,1)^n$ converges
	asymptotically to a local minimizer of $E$ with probability one.
\begin{thm}\longthmtitle{Convergence of $\binpac$}\label{thm:cent-cvg}
Given an initial condition $x(0)\in(0,1)^n$, the
    trajectory $x(t)$ under $\binpac$ converges asymptotically to a
    critical point $x^\star$ of $E$. In addition, under
  Assumption~\ref{ass:init}, on the random choice of initial
  conditions, and Assumption~\ref{ass:t-tau}, on the random choice of
  $T,\tau$, the probability that $x(0)$ is in the set
  $\underset{\hat{x}}{\cup} \W^s (\hat{x})$, where $\hat{x}$ is a
  saddle-point or local maximum of $E$, is zero.
\end{thm}

We now define a Deterministic Annealing (DA) variant inspired
by~\cite{KR:98} to augment the $\binpac$ dynamics and provide a method
for gradually learning a justifiably good feasible point of $\P 1$. In~\cite{KR:98}, the author justifies the deterministic online tuning of a temperature parameter in the context of data clustering and shows that this avoids poor local optima by more thoroughly exploring the state space. Similarly, we aim to learn a sufficiently good solution trajectory by allowing the dynamics to explore the interior of the unit hypercube in the early stages of the algorithm, and then to force the trajectory outward to a feasible binary solution by gradually adjusting $T$ or $\tau$ online. 

Consider either reducing the temperature $T$ or increasing the time
constant $\tau$ during the execution of $\binpac$. 
This reduces the terms in $E$ which promote
convexity, particularly near the boundaries of the unit hypercube. As
$T,\tau$ are adjusted, for $a\prec -\gamma\Vert p\Vert^2$, the domain
of $E$ becomes gradually more concave away from the corners of the
unit hypercube. Thus, starting with $T/\tau$ sufficiently large, the
early stages of the algorithm promote exploration of the interior of
the state space. As $T/\tau$ is reduced, the equilibria of $E$ are
pushed closer to (and eventually converge to) the feasible points of
$\P 1$. The update policy we propose is described formally in
Algorithm~\ref{alg:da}, and we further explore its performance in simulation. 

\begin{algorithm}
	\caption{Determinisitc Annealing}\label{alg:da}
	\begin{algorithmic}[1]
		\Procedure{Det-Anneal}{$\beta > 1,T_0,\tau_0,t_d$}
		\State Initialize $x(0)$
		\State $T\gets T_0, \tau\gets \tau_0$
		\While{\texttt{true}}
		\State Implement $\binpac$ for $t_d$ seconds\label{line:da-condition}
		\State $\tau \gets \beta \tau \text{\quad or\quad} T \gets (1/\beta) T$
		\EndWhile
		\EndProcedure
	\end{algorithmic}
\end{algorithm}

Note that Algorithm~\ref{alg:da} leads to a hybrid
  dynamic system with discrete jumps in an enlarged state $\phi =
  (x,T,\tau)$, which can cast some doubt on basic existence and
  uniqueness of solutions. We refer the reader to Propositions 2.10
  and 2.11 of~\cite{RG-RS-AT:09} to justify existence and uniqueness
  of solutions in the case of $t_d > 0$ fixed.

\begin{cor}\longthmtitle{Convergence to Feasible Points}
  Under Assumptions~\ref{ass:init}-\ref{ass:t-tau} and $a \prec -\gamma\| p\|^2$, the $\binpac$
  dynamics augmented with Algorithm~\ref{alg:da} converge
  asymptotically to feasible points of $\P 1$.
\end{cor}

The result of the Corollary is quickly verified by inspecting the
terms of $H(x)$. 
The function $E$ 
is smooth, strictly concave near $x=(0.5)\ones_n$ for small $T/\tau$ due to the design of $a_i$, and becomes strictly convex as the elements of $x$ approach
$0$ or $1$, corresponding
to isolated local minima of $E$, due to the $T/\tau$ term dominating $H(x)$. As the quantity
$T/\tau$ is reduced under Algorithm~\ref{alg:da}, these local minima
are shifted asymptotically closer to corners of the unit hypercube, i.e.
feasible points of $\P 1$.

\section{Distributed Hopfield Neural Network}\label{sec:dist-hop}

With the framework of the previous section we formulate a problem $\P
2$ which is closely related to $\P 1$, but for which the global
penalty term can be encoded by means of an auxiliary decision
variable. This formulation leads to the Distributed $\binpromd$, or $\binpad$, which we rigorously analyze for its
convergence properties.

It is clear from the PT-inverse operation and $W$ being nonsparse that
$\binpac$ is indeed centralized. In this section, we design a
distributed algorithm in which each agent $i$ must only know
$p_j, j\in\N_i$ and the value of an auxiliary variable $y_j, j\in\N_i
\cup\N_i^2$, i.e. it must have communication with its two-hop neighbor
set. If two-hop communications are not directly available, the
algorithm can be implemented with two communication rounds per
algorithm step. We provide comments on a one-hop algorithm in
Remark~\ref{rem:one-hop}.

\begin{assump}\longthmtitle{Graph Properties and Connectivity}\label{ass:graph-conn}
  The graph $\mathcal{G} = (\N,\mathcal{E})$ is undirected and
  connected; that is, a path exists between any two pair of nodes and,
  equivalently, its associated Laplacian matrix $L = L^\top$ has rank
  $n-1$.
\end{assump}
Now consider the $n$ linear equations $(p_i x_i)_i + Ly =
(\subscr{P}{r}/n)\ones_n.$ Notice that, by multiplying from the left
by $\ones_n^\top$ and applying $\ones_n^\top L = \zeros_n^\top$, we
recover $p^\top x = \subscr{P}{r}$. Thus, by augmenting the state with
an additional variable
$y\in\real^n$, 
we can impose a distributed penalty term. We now formally state the
distributed reformulation of~$\P 1$:
\begin{equation*} 
  \P 2: \
  \underset{x\in\{0,1\}^n,y\in\real^n}{\text{min}} \
  \tilde{f}(x,y) = \sum_i^n f_i(x_i) + \dfrac{\gamma}{2}\sigma^\top \sigma,  
\end{equation*}
where the costs $f_i$ again satisfy $f_i(1) - f_i(0) = c_i$ and we
have defined $\sigma = (p_i x_i)_i + Ly - (\subscr{P}{r}/n) \ones_n$ for
notational simplicity. Before proceeding, we provide some context on
the relationship between $\P 1$ and $\P 2$.

\begin{lem}\longthmtitle{Equivalence of P1 and P2}\label{lem:equiv}
  Let Assumption~\ref{ass:graph-conn}, on graph connectivity, hold, and
  let $(x^\star,y^\star)$ be a solution to $\P 2$. Then, $x^\star$ is
  a solution to $\P 1$ and $f(x^\star) = \tilde{f}(x^\star,y^\star)$.
\end{lem}

To define $\binpad$, we augment the centralized $\binpac$
with gradient-descent dynamics in $y$ on a newly obtained energy function $\widetilde{E}$ of $\P 2$. Define $\widetilde{E}$ as
\begin{equation}\label{eq:energyd}
\widetilde{E}(x,y) = \tilde{f}(x,y) + \dfrac{1}{\tau} \sum_i \int_0^{x_i} g^{-1}(\nu) d\nu.
\end{equation}

In Section~\ref{sec:all-hop}, we obtained a matrix $W$ which was
nonsparse. Define $\widetilde{W}, \tilde{v}$ for $\widetilde{E}$ via $\tilde{f}$ as $\widetilde{W} = -\diag{a + \gamma (p_i^2)_i},
\tilde{v} = (a_i b_i)_i + \gamma\diag{p}\left((\subscr{P}{r}/n)\ones_n
- Ly\right).$
Compute the Hessian of $\widetilde{E}$ with respect to only $x$ as $\widetilde{H}(x) = \nabla_{xx}\widetilde{E}(x,y) = -\widetilde{W} + (T/\tau)\diag{1/x-(x_i^2)_i}$.
Since $\widetilde{H}(x)$ is diagonal, the $\supscr{ii}{th}$ element of
the PT-inverse of $\widetilde{H}(x)$ can be computed locally by each
agent $i$ as:
\begin{equation*}
  (\vert \widetilde{H}(x)\vert_m)^{-1}_{ii} = \begin{cases}
    \vert \widetilde{H}(x)_{ii} \vert^{-1}, 
    &  \vert\widetilde{H}(x)_{ii}\vert \geq m, \\
    1/m, & \text{o.w.}
\end{cases}
\end{equation*}
where $\widetilde{H}(x)_{ii} = a_i + \gamma p_i^2 + T/\tau (x_i -
x_i^2)^{-1}$. 
The $\binpad$ dynamics, which are PT-Newton descent in
$x$ and gradient descent in $y$ on $\widetilde{E}$, are then stated
as:
\begin{equation}\label{eq:binpad}
	\begin{aligned}
          \dot{x} &=(\vert \widetilde{H}(x)\vert_m)^{-1} \diag{\frac{(x_i-x_i^2)_i}{T}} \\
          &\qquad \qquad \qquad \left(\widetilde{W}x + \dfrac{T}{\tau}\log\left(1/x_i-1\right)_i + \tilde{v}\right), \\
          \dot{y} &= -\alpha\gamma
          L\left(
          	(p_i x_i)_i
            + Ly\right),
	\end{aligned}
\end{equation}
Due to the new matrices $\widetilde{W},\tilde{v}$ and the sparsity of
$L$, $\dot{x}$ can be computed with one-hop information and $\dot{y}$
with two-hop information (note the $L^2$ term);
thus,~\eqref{eq:binpad} defines a distributed algorithm. Additionally,
recalling the discussion on parameter design, the problem data $a$ and
$b$ can now be locally designed.

Before proceeding, we establish a property of the domain of $y$ and some distributed extensions of Lemmas~\ref{lem:finite-eq} and~\ref{lem:fwd-inv-c}.
\begin{lem}\longthmtitle{Domain of Auxiliary Variable}\label{lem:y-domain}
	Given an initial condition $y(0)$ with $\ones_n^\top y(0) = \kappa$, the trjaectory $y(t)$ is contained in the set
	\begin{equation}\label{eq:ydom}
	\Y = \setdef{\omega + (\kappa/n)\ones_n}{\ones_n^\top \omega = 0}.
	\end{equation}
\end{lem}

\begin{lem}\longthmtitle{Closed Form Auxiliary Solution}\label{lem:y-soln}
	For an arbitrary fixed $x\in[0,1]^n$, the unique minimizer $y^\star$ contained in $\Y$ of both $\tilde{f}$ and $\widetilde{E}$ is given by
	\begin{equation}\label{eq:ystar}
	y^\star = -L^\dagger \left(
	p_i \tilde{x}_i\right)_i + \frac{\kappa}{n}\ones_n.
	\end{equation}
	This is also the unique equilibrium of~\eqref{eq:binpad} in $\Y$.
\end{lem}

\begin{lem}\longthmtitle{Finite Equilibria (Distributed)}~\label{lem:finite-eq-dist}
	Let $\widetilde{\X}\times\widetilde{\Y}$ be the set of equilibria of~\eqref{eq:binpad}
	satisfying $(\dot{x},\dot{y}) = 0$ on $(x,y)\in[0,1]^n \times \Y$. The set $\widetilde{\X}\times\widetilde{\Y}$ is finite.
\end{lem}

We now extend the results of Theorem~\ref{thm:cent-cvg} to the
distributed case of solving $\P 2$ via $\binpad$. 
We have the following theorem on the trajectories of
$(x(t),y(t))$
under~\eqref{eq:binpad}, which can be interpretted as establishing convergence to a local minimizer with probability one.

\begin{thm}\longthmtitle{Convergence of $\binpad$}\label{thm:dist-cvg}
  Given an initial condition $(x(0),y(0))\in(0,1)^n\times \real^n$,
 the trajectory $(x(t),y(t))$ under $\binpad$
  converges asymptotically to a critical point $(x^\star,y^\star)$ of
  $\widetilde{E}$. In addition, under Assumption~\ref{ass:init}, on
  the random choice of initial condition $x(0)$, and
  Assumption~\ref{ass:t-tau}, on the random choice of $T,\tau$,
  the probability that $(x(0),y(0))$ is in the set
    $\underset{\hat{x},\hat{y}}{\cup} \W^s (\hat{x},\hat{y})$, where
    $(\hat{x},\hat{y})$ is a saddle-point or local maximum of
    $\widetilde{E}$, is zero. Lastly, all local minima
  $(x^\star,y^\star)$ of $\widetilde{E}$ are globally optimal in $y$:
  $\widetilde{E}(x^\star,y)\geq \widetilde{E}(x^\star,y^\star),
  \forall y\in\real^n$.
\end{thm}

\begin{lem}\longthmtitle{Forward Invariance of the Open Hypercube (Distributed)}\label{lem:fwd-inv-d}
	The set $(0,1)^n \times \Y$ is a forward-invariant set under the
	$\binpad$ dynamics~\eqref{eq:binpad}.
\end{lem}

\begin{rem}\longthmtitle{One-Hop Distributed
    Algorithm}\label{rem:one-hop} \rm The proposed distributed algorithm requires
  two-hop neighbor information, which may be intractable in some
  settings. The source of the two-hop term stems from the quadratic
  $\gamma$ penalty term. However, it is possible to define a one-hop
  distributed algorithm via a Lagrangian-relaxation route.
	
  Consider posing $\P 2$ with the $\gamma$ term instead as a linear
  constraint: $\sqrt{\gamma /2} ((p_i x_i)_i + Ly) = \sqrt{\gamma /2}
  (\subscr{P}{r}/n)\ones_n$. Applying Lagrangian
    relaxation to this problem introduces a Lagrange multiplier on
  the linear terms, and from there it would be appropriate to define a
  saddle-point-like algorithm along the lines
  of~\cite{AC-EM-SHL-JC:18-tac} in which gradient-ascent in the dual
  variable is performed. This changes the nature of the penalty from
  squared to linear, so the underlying optimization model is different
  in that sense, but it follows that this approach could be
  implemented with one-hop information.
  
  We note that, in some distributed contexts, penalty terms or constraints can be imposed via $\sqrt{L}$ which then appears as $L$ in the associated squared terms of the dynamics (in place of $L^2$). However, the linear $L$ also appears in our algorithm, and substituting $\sqrt{L}$ would not inherit the sparsity of the communication graph. Therefore we leave the design of a fully one-hop mixed first-order/second-order algorithm as an open problem. 
\end{rem}

\section{Simulations}\label{sec:sims}

Our simulation study is split in to two parts; the first focuses on 
numerical comparisons related to runtime and solution quality, and
the second is a 2D visualization of the trajectories of the
Distributed Annealing (DA) variants for both the centralized and
distributed NNN methods.

\subsection{Runtime and Solution Quality Comparison}\label{ssec:runtime-soln}

In this section, we compare to a greedy method stated as
Algorithm~\ref{alg:greedy} and a semidefinite programming (SDP)
relaxation method stated as Algorithm~\ref{alg:sdp}. In short, the
greedy method initializes the state as $x = \zeros_n$ and iteratively
sets the element $x_i$ to one which decreases the cost function the
most. This is repeated until no element remains for which the updated
state has lower cost than the current state. For the SDP method, a
convex SDP is obtained as the relaxation of $\P 1$, see
e.g.~\cite{LV-SB:96}. We use the shorthand $\texttt{SDPrlx}(\bullet)$
to indicate this in the statement of Algorithm~\ref{alg:sdp}. This SDP is solved
using CVX software in MATLAB~\cite{website:cvx} and a lowest-cost
partition is computed to construct a feasible solution. For the
sake of convenience in stating both algorithms, we have defined $f':
2^n\rightarrow\real$ to be the set function equivalent of $f$,
i.e. the cost of $\P 1$. That is, $f'(\SSS) = f(x)$, where $i\in\SSS$
indicates $x_i = 1$ and $i\notin\SSS$ indicates $x_i = 0$. Finally, we additionally compare to a brute force method which we have manually programmed as an exhaustive search over the entire (finite) feasibility set.

\begin{algorithm}
	\caption{Greedy Method}\label{alg:greedy}
	\begin{algorithmic}[1]
		\Procedure{Greedy}{$f'$}
		\State $\SSS \gets \emptyset$
		\State $\texttt{done} \gets \texttt{false}$
		\While{$\texttt{done} = \texttt{false}$}
		\State $i^\star \gets \underset{i\notin\SSS}{\argmin} \ f'(\SSS \cup \{i\})$
		\If{$f'(\SSS\cup \{i^\star\}) < f'(\SSS)$}
		\State $\SSS \gets \SSS \cup \{i^\star\}$
		\Else
		\State $\texttt{done} \gets \texttt{true}$
		\EndIf
		\EndWhile
		\State $x_i \gets \begin{cases}
		0, & i\notin \SSS, \\
		1, & i\in \SSS.
		\end{cases}$
		\State \textbf{return} $x$
		\EndProcedure
	\end{algorithmic}
\end{algorithm}

\begin{algorithm}
	\caption{SDP Relaxation Method}\label{alg:sdp}
	\begin{algorithmic}[1]
		\Procedure{SDP}{$f'$}
		\State $\P_{\text{SDP}} \gets \texttt{SDPrlx}(\P 1)$
		\State $x^\star \gets \underset{x}{\argmin} \P_{\text{SDP}}$
		\State $\SSS \gets \emptyset$
		\State $\texttt{done} \gets \texttt{false}$
		\While{$\texttt{done} = \texttt{false}$}
		\State $i^\star \gets \underset{i\notin\SSS}{\argmax} \ x_i$
		\If{$f'(\SSS\cup \{i^\star\}) < f'(\SSS)$}
		\State $\SSS \gets \SSS \cup \{i^\star\}$
		\Else
		\State $\texttt{done} \gets \texttt{true}$
		\EndIf
		\EndWhile
		\State $x_i \gets \begin{cases}
		0, & i\notin \SSS, \\
		1, & i\in \SSS.
		\end{cases}$
		\State \textbf{return} $x$
		\EndProcedure
	\end{algorithmic}
\end{algorithm}

%

\begin{figure}[h]
	\centering
	\includegraphics[scale = 0.4]{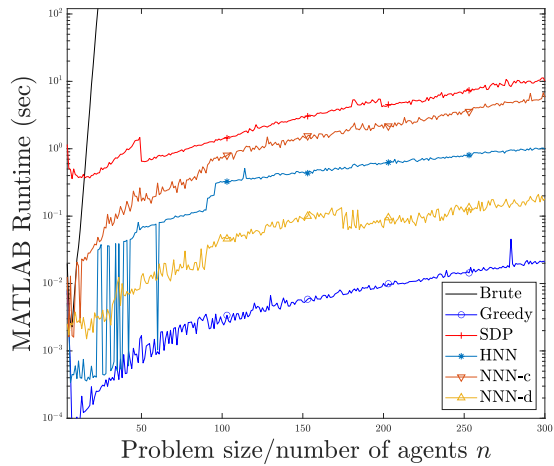}
	\caption{Runtime of each method for increasing problem sizes.}
	\label{fig:runtime}
\end{figure}

In Figure~\ref{fig:runtime} 
we plot the runtime in
MATLAB on a 3.5GHz Intel Xeon E3-1245 processor over increasing problem size
$n$ for each of six methods: a brute force search, the aforementioned
greedy and SDP methods, the HNN first proposed in~\cite{JH-DT:85}
(i.e. the gradient-like version of $\binpac$), and the $\binpac$ and
$\binpad$ methods we developed in Sections~\ref{sec:all-hop}
and~\ref{sec:dist-hop}. The first obvious observation to make is that
the runtime of brute force method increases at a steep exponential
rate with increasing $n$ and exceeds 120 seconds at $n=22$, making it
intractable for even medium sized problems. Next, we note that there
are some spikes associated with the HNN method around $n=25$
to $n=40$. These are reproducible, and we suspect that this is due to
the emergence of saddle-points and increasing likelihood of
encountering these along the trajectory as $n$ increases. This is a well-documented problem observed in literature, see e.g.~\cite{YD-RP-CG-KC-SG-YB:14}, and we also confirm it empirically in this setting by observing that share of iterations for which the Hessian is indefinite (as opposed to positive definite) tends to grow as $n$ increases. We also note
that $\binpac$ scales relatively poorly, which can be attributed to a
matrix eigendecomposition being performed at each discretized
iteration of the continuous-time algorithm. For $\binpad$, the matrix
being eigendecomposed is diagonal, which makes it a trivial operation
and allows $\binpad$ to scale well. We note that the SDP method scales the worst amongst the non brute-force methods. Unsurprisingly, the
greedy method remains the fastest at large scale, although recall that
the motivation of developing our method is for it to be distributed
and that a greedy approach can not be distributed due to the global penalty term.

As for algorithm performance as it pertains to the cost of the
obtained solution, we fix $n=50$ and additionally include DA variants of both $\binpac$ and $\binpad$. We also omit the brute force method due
to intractability. For the sake of comparison, we compute a
performance metric $Q$ and provide it for each method in
Table~\ref{table:perf}.
The metric $Q$ is computed
as follows: for each trial, sort the methods by solution cost. Assign
a value of 6 for the best method, 5 for the second-best, and so on,
down to the seventh-best (worst) receiving zero. Add up these scores
for all 100 trials, and then normalize by a factor of 600 (the maximum
possible score) to obtain $Q$. Note that $Q$ does not account for runtime in any way.

It should be unsurprising that the
tried-and-true centralized greedy and SDP methods perform the
best. However, we note that they were beaten by our methods in a
significant number of trials, which can be seen by noting that a $Q$ score for
two methods which perform best or second-best in all trials would sum
to $1100/600 = 1.83$, while $Q(\text{greedy})+Q(\text{SDP})=1.75$, or a
cumulative pre-scaled score of $1050$, indicating that our methods
outperformed these methods in net 50 ``placement spots" over the 100
trials. In general, we find that the DA version of the NNN algorithms
obtains better solutions than the non-DA version, confirming the
benefit of this approach. We also find that 
$\binpad$ generally outperforms $\binpac$.
It's possible that an initially ``selfish" trajectory in $x$ is beneficial, which would
neglect the global penalty until $y$ adequately converges, although this is speculative. Lastly,
we note that the HNN method never performs better than worst, which we
attribute to the steepest-descent nature of gradient algorithms which
do not use curveature information of the energy function. It might be
possible that the stopping criterion forces HNN to terminate near
saddle-points, although we do not suspect this since we observe the
Hessian is positive-definite in the majority of termination
instances.

As for parameter selection, we find that choosing $m \ll 1$ is generally
best, since $m\geq 1$ would always produce a PT-inverse Hessian with
eigenvalues contained in $(0,1]$. This effectively scales down
$\dot{x}$ in the eigenspace associated with Hessian eigenvalue
magnitudes greater than $1$, but does not correspondingly scale up
$\dot{x}$ in the complementary eigenspace associated with small
eigenvalues. Additionally, choosing $T/\tau$ greater than
$1$ in the fixed case 
tended to be
effective. This may be related to selecting $a_i < -\gamma\Vert
p\Vert^2 - 4T/\tau$ to guarantee anti-stability from $(0.5)\ones_n$,
and would explain why a high $T_0/\tau_0$ that decreases in the
DA learning variant performs so well. In general, for the DA learning variant, we recommend choosing $T_0,\tau_0$ so that $T_0/\tau_0 \gg 1$ and also $\beta>1$ sufficiently large so that $T/\tau \ll 1$ by algorithm termination, which gives rise to a robust exploration/exploitation tradeoff. Finally, all $\alpha\approx 1$
seem to behave roughly the same, with only $\alpha \ll 1$ and
$\alpha \gg 1$ behaving poorly (the former leading to slow convergence
in $y$ and ``selfish" behavior in $x$, and the latter being
destabilizing in the discretization of $\dot{y}$).


\begin{table}[]
	\centering
	\caption{Comparison of performance metric $Q$ for $100$ randomized trials with $n=50$.} 
	\label{table:perf}
	\scalebox{1}{
		\begin{tabular}{|c|c|c|c|c|}
			\hline
			Method & $Q$ \\ \hline \hline
			$\binpac$ & 0.2161 \\ \hline
			$\binpac$-DA & 0.2891 \\ \hline
			$\binpad$ & 0.5443 \\ \hline
			$\binpad$-DA & 0.7005 \\ \hline
			HNN & 0 \\ \hline
			Greedy & 0.8411 \\ \hline
			SDP & 0.9089 \\ \hline
	\end{tabular}}
\end{table}

\begin{table}[]
	\centering
	\caption{Problem data and parameter choices (where relevant) for performance comparison. Problem data $p_i, c_i$ is generated randomly from given distributions for each of $100$ trials.}
	\label{table:data}
	\scalebox{1}{
		\begin{tabular}{|c|c|c|c|c|}
			\hline
			Data or parameter & Value \\ \hline \hline
			$n$ & $50$ \\ \hline
			$p_i$ & $\U[1,50]$ \\ \hline
			$c_i$ & $p_i^e$, $e\sim\U[2,3]$ \\ \hline
			$\subscr{P}{r}$ & $1500$ \\ \hline
			$\gamma$ & $1$ \\ \hline
			$T_0$ & $1$ \\ \hline
			$\tau_0$ & $0.1$ \\ \hline
			$m$ & $0.1$ \\ \hline
			$\alpha$ & $1$ \\ \hline
			Learning steps & $10$ \\ \hline
			$\beta$ & $1.4$ \\ \hline
			$n$ & $50$ \\ \hline
	\end{tabular}}
\end{table}

\subsection{Learning Steps and 2-D Trajectories}\label{ssec:traj-da}

\begin{figure}[h]
	\centering
	\subfloat[]{{\includegraphics[scale=.55]{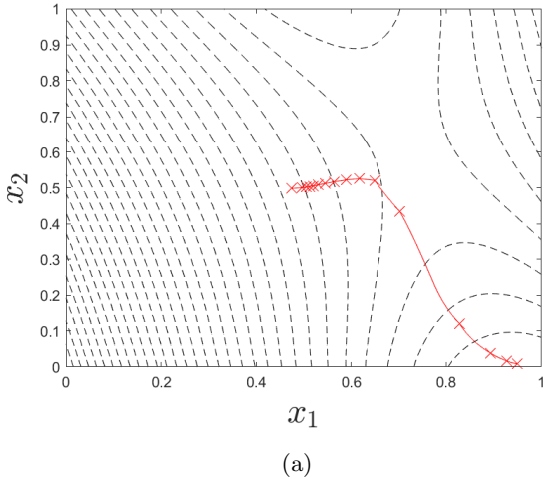} }}\label{fig:centnewt} \\
	\subfloat[]{{\includegraphics[scale=.55]{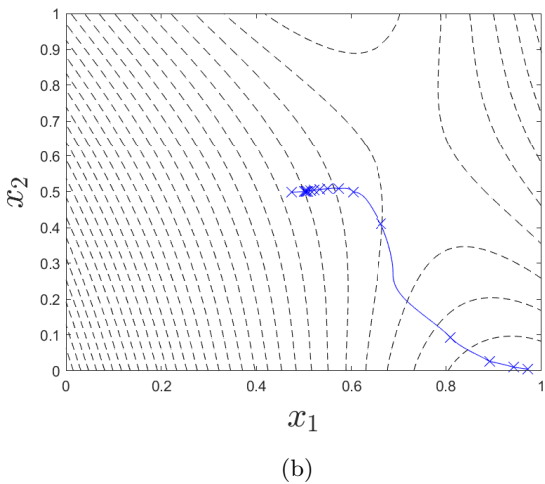} }}\label{fig:distnewt}
	\caption{Centralized $\binpac$ (a) and distributed $\binpad$ (b) trajectories in 2D with $15$ learning steps. Stable equilibrium points between learning steps indicated by $\times$, contours of $E$ and $\widetilde{E}$ in final step indicated by dashed lines.}
	\label{fig:traj}
\end{figure}

Next, for the sake of understanding how the learning rate $T/\tau$ affects the trajectories of the solutions, we have provided Figure~\ref{fig:traj} which plots the 2-D trajectories of $\binpac$ and $\binpad$ with $T/\tau$ being gradually reduced over 15 learning steps. The contours of the energy function for the final step are also plotted. The problem data and choice for $a$ is:
\begin{equation*}
\begin{aligned}
c &= (2,1)^\top, \quad p = (3,1)^\top, \quad \subscr{P}{r} = 2.8, \quad \gamma = 4, \\
a &:= -(10,10)^\top. 
\end{aligned}
\end{equation*}
Note that, in each case, the trajectory approaches the optimal solution $x^\star = (1,0)^\top$. However, it is worth noting that a steep saddle point occurs around $x=(0.75,0.6)^\top$. Intuitively, this corresponds to a high risk of the trajectory veering away from the optimal solution had the DA not been implemented. With the opportunity to gradually learn the curveature of the energy function, as shown by stabilization to successive equilibria marked by $\times$, each algorithm is given the opportunity to richly explore the state space before stabilizing to the optimal solution $(1,0)^\top$. Further studying the learning-rate $T/\tau$ and a more complete analysis of Algorithm~\ref{alg:da} and the parameter $\beta$ are subjects of future work.

\section{Conclusion}

This paper posed an optimal generator dispatch problem for settings in
which the agents are generators with binary controls. We first showed that the centralized problem is amenable to
solution via a Centralized $\binpromc$ approach and
proved convergence to a local minimizer with probability one
under light assumptions. Next, we developed an
approach to make the dynamics computable in a distributed
setting in which agents exchange messages with their two-hop
neighbors in a communication graph. The methods scale and perform well compared to
standard greedy and SDP-relaxation approaches, and the latter method enjoys the qualities of a distributed algorithm, 
unlike previous approaches. Future research
directions include application of the methods to a broader class of problems which may include additional cost terms or constraints and a deeper analysis of the Deterministic Annealing variant as it pertains to the online adjustment of the learning-rate $T/\tau$.

\section*{Appendix}
\subsection*{Proofs of Lemmas, Theorems, and Propositions}

\textbf{Proof of Lemma~\ref{lem:finite-eq}:}
	First consider only $\X\cap (0,1)^n$. Note that $(\vert H(x)\vert_m
	)^{-1} \succ 0$ (by construction) and $\diag{(x_i-x_i^2)_i/T}\succ 0$ on
	$x\in(0,1)^n$, so we focus on
	\begin{equation}\label{eq:grad-zero}
	Wx + \frac{T}{\tau}\log\left(1/x_i-1\right)_i + v = \zeros_n.
	\end{equation}
	Examining the above expression elementwise, it is nonconstant,
	continuous, and its derivative changes sign only a finite
	number of times. Therefore, the total number of zeros on
	$(0,1)^n$ must be finite.
	
	Now consider the $\supscr{i}{th}$ element
	of~\eqref{eq:grad-zero} for $x_j\rightarrow 0$ or $1$ for all
	$j$ in an arbitrary permutation of
	$\until{n}\setminus{\{i\}}$. Since the number of these
	permutations is finite, and each permutation still gives rise
	to a finite number of solutions to~\eqref{eq:grad-zero} in the
	$\supscr{i}{th}$ component, it follows that $\X$ is finite.
\QED

\textbf{Proof of Lemma~\ref{lem:fwd-inv-c}:}
	Consider again the terms of $\dot{x}$ elementwise. There are two cases to
	consider for evaluating $x_i$: $x_i = \varepsilon$ and $x_i =
	1-\varepsilon$ for some $0< \varepsilon \ll 1$ sufficiently small
	such that the terms of $(\vert H(x)\vert_m)^{-1}$ are still
	dominated by $(1/x_i-x_i^2)$ and the $Wx + v$ are still dominated by
	the $\log$ term. Then, consider the expression
	\begin{equation}\label{eq:log-express}
	\log\left( 1/x_i - 1\right)(x_i-x_i^2)^2.
	\end{equation}
	For $x_i = \varepsilon \approx 0$,~\eqref{eq:log-express} evaluates
	to a small positive value, and for $x_i = 1-\varepsilon \approx
	1$,~\eqref{eq:log-express} evaluates to a small negative value. We
	have argued that these are the dominating terms regardless of values
	of the remaining components of $x$, and so we conclude that
	$x_i\in\{0,1\}$ are componentwise anti-stable and that elements of
	$x$ will never approach $0$ or $1$. Thus, the open hypercube is
	forward invariant.
\QED

\textbf{Proof of Theorem~\ref{thm:cent-cvg}:}
	Let $\X$ be the set of all critical points of $E$. We first
	establish that $E$ decreases along the trajectories of $\binpac$ and
	that $x(t)$ converges asymptotically to $\X$. Differentiating $E$ in
	time, we obtain:
	\begin{equation}\label{eq:dE-dt}
	\begin{aligned}
	\dfrac{dE}{dt} &=
	\dot{x}^\top \nabla_x E(x) = \dot{x}^\top\left(-Wx -v + g^{-1}(x)/\tau\right) \\
	&= -\dot{x}^\top \diag{\frac{T}{(x_i-x_i^2)_i}}\vert H(x)\vert_m \dot{x} <0, \\
	& \qquad \qquad \qquad \text{for} \ \dot{x} \neq 0, \
	x\in(0,1)^n.
	\end{aligned}
	\end{equation}
	Recall that $x(t)\in(0,1)^n$ for all $t\geq 0$ due to
	Lemma~\ref{lem:fwd-inv-c}. From~\eqref{eq:binpac} and the discussion
	that followed on equilibria, $\dot{x} = 0$ implies $\nabla_x E(x) =
	0$ due to $(\vert H(x) \vert_m)^{-1} \succ 0$ and $\diag{(x_i-x_i^2)_i/T}
	\succ 0$ on $x\in(0,1)^n$. The domain of $E$ is the compact set
	$[0,1]^n$ (per the definition of the integral terms), and $E$ is
	continuous and bounded from below on this domain, so at least one
	critical point exists. Combining this basic fact
	with~\eqref{eq:dE-dt} shows that the $\binpac$ dynamics
	monotonically decrease $E$ until reaching a critical point. More
	formally, applying the LaSalle Invariance Principle\cite{HKK:02}
	tells us that the trajectories converge to the largest invariant set
	contained in the set $dE/dt = 0$. This set is $\X$, which is finite
	per Lemma~\ref{lem:finite-eq}. In this case, the LaSalle Invariance Principle
	additionally establishes that we converge to a single
	$x^\star\in\X$.

	The proof of the second statement of the theorem relies on an
	application of the Stable Manifold Theorem (see~\cite{JG-PH:83}) as well as Lemma~\ref{lem:finite-eq}. Let $\dot{x} = \varphi_{T,\tau}(x)$
	for a particular $T,\tau$. 
	We aim to show that 
	$\Pp[\cup_{\hat{x}} \setdef{\W_s (\hat{x})}{\hat{x} \text{ is a saddle
			or local maximum}} ] = 0$
	under Assumptions~\ref{ass:init}-\ref{ass:t-tau}. It is
	sufficient to show that, for each critical point $x^\star$
	such that $\varphi_{T,\tau}(x^\star) = 0$, and almost all
	$T,\tau$, $D\varphi_{T,\tau} (x^\star)$ is full rank and
	its eigenvalues have non-zero real
	parts.  
	The reason for this argument is the following: let $x^\star$ be a
	critical point with $D\varphi_{T,\tau} (x^\star)$ full rank and
	eigenvalues with non-zero real parts.  If the eigenvalues do
	\emph{not} all have positive real parts, then some have negative
	real parts, which indicates that $x^\star$ is a saddle or local
	maximum of $E$. These negative real-part eigenvalues induce an
	unstable manifold of dimension $n_u \geq 1$. As such, the globally
	stable set $\W_s (x^\star)$ is a manifold with dimension $n - n_u
	< n$, and $\Pp\left[x(0)\in\W_s(x^\star) \right]=0$ per
	Assumption~\ref{ass:init}.
	
	To argue this case, define $h:(0,1)^n\times \real\times\real
	\rightarrow \real$ as
	\begin{equation*}
	h(x,T,\tau) = \det{D\varphi_{T,\tau}(x)}.
	\end{equation*}
	We now leverage Assumption~\ref{ass:t-tau} and~\cite{BM:15} to claim
	first that $\bar{\Pp}\left[h(x^\star,T,\tau) = 0\right] = 0$
	for each $x^\star\in\X$, i.e. $D\varphi_{T,\tau} (x^\star)$ is full
	rank for each $x^\star$ with probability one w.r.t.~$\bar{\Pp}$.  We
	first address the points $x$ for which the function $h$ is
	discontinuous. Define $\hat{\X}$ as the set of $x$ for which the
	truncation of the eigenvalues of $H(x)$ becomes active, i.e. the
	discontinuous points of $h$. Although we do not write it as such,
	note that $H$ is implicitly a function of $T,\tau$ and that the
	eigenvalues of $H$ can be expressed as nonconstant real-analytic
	functions of $T,\tau$. Considering this fact and an arbitrary $x$,
	the set of $T,\tau$ which give $x\in\hat{\X}$ has measure zero with
	respect to $\real^2$\cite{BM:15}. Thus, for particular $T,\tau$, $h$
	is $C^{\infty}$ almost everywhere. Applying once more the argument
	in~\cite{BM:15} and Assumption~\ref{ass:t-tau} with the fact that
	$h$ is a nonconstant real analytic function of $T,\tau$ we have that
	\begin{align*}
	\bar{\Pp}\,[\mathcal{T}(\hat{x}) \triangleq  \{&(T,\tau) \,| \,
	h(\hat{x}, T,\tau) = 0 \}]  = 0, \quad \forall\, \hat{x}\notin \hat{\X}.
	\end{align*} 
	Now consider the set of critical points as an explicit function of
	$T,\tau$ and write this set as $\X(T,\tau)$. 
	Recalling Lemma~\ref{lem:finite-eq}, the set of $\hat{x}$ that we
	are interested in reduces to a finite set of critical points
	$x^\star\in \X(T,\tau)$. Thus, we can conclude that
	$\bar{\Pp}(\cup_{x^\star \in \X(T,\tau)}\,\mathcal{T}(x^\star))
	\le \sum_{x^\star \in \X(T,\tau)} \bar{\Pp}(\mathcal{T}(x^\star))=
	0$.
	
	There is an additional case which must be considered, which is that
	$h(x^\star ,T,\tau)\neq 0$, but some eigenvalues of
	$D\varphi_{T,\tau}(x^\star)$ are purely imaginary and induce stable
	center manifolds, which could accommodate the case of a globally
	stable set which is an $n$-dimensional manifold (i.e. the ``degenerate saddle" case). We consider the
	function $h$ mostly out of convenience, but the argument can be
	extended to a function $\mathbf{h}: (0,1)^{n}\times \real \times
	\real \rightarrow \complex^n$ which is a map to the roots of the
	characteristic equation of $D\varphi_{T,\tau}(x)$. We are concerned
	that each element of $\mathbf{h}(x,T,\tau )$ should have a nonzero
	real part almost everywhere. To extend the previous case to this,
	consider the identification $\complex \equiv \real^2$ and compose
	$\mathbf{h}$ with the nonconstant real analytic function $\zeta(w,z)
	= w$, for which the zero set is $w \equiv 0$, corresponding to the
	imaginary axis in our identification. From this, we obtain a
	nonconstant real-analytic as before whose zero set is the imaginary
	axis. Applying the argument in~\cite{BM:15} in a
	similar way as above, $\mathbf{h}(x,T, \tau )$ has nonzero real
	parts for almost all $(T,\tau )$ for each
	$x$. Therefore, the probability of a particular saddle point or local maximum
	$x^\star$ having a nonempty stable center manifold is zero for
	arbitrary $x(0)$ satisfying Assumption~\ref{ass:init} and $T,\tau$ satisfying
	Assumption~\ref{ass:t-tau}.
	\QED

\textbf{Proof of Lemma~\ref{lem:equiv}:}
	The equivalence stems from the global term and the flexibility in
	the unconstrained $y$ variable. Notice
	\begin{equation*}
	\begin{aligned}
	\dfrac{\gamma}{2} \sigma^\top \sigma & = \dfrac{\gamma}{2}
	\sigma^\top (I_n -
	\ones_n\ones_n^\top/n)\sigma + \dfrac{\gamma}{2} \sigma (\ones_n\ones_n^\top/n)\sigma \\
	&= \dfrac{\gamma}{2} \sigma^\top (I_n -
	\ones_n\ones_n^\top/n)\sigma + \dfrac{\gamma}{2}(p^\top x -
	\subscr{P}{r})^2.
	\end{aligned}
	\end{equation*}
	We have recovered the original global term of $\P 1$ in the
	bottom line, so now we deal with the remaining term. The
	matrix $I_n - \ones_n \ones_n^\top/n \succeq 0$ has $\image
	I_n - \ones_n \ones_n^\top/n = \spn \{\ones_n\}^\perp = \image
	L$, given that $L$ is connected. Thus, because $y$ is
	unconstrained and does not enter the cost anywhere else, we
	can compute the set of possible minimizers of $\tilde{f}$
	in closed form with respect to any $x$ as
	\begin{equation*}
	\begin{aligned}
	y^\star &\in \setdef{-L^\dagger \big((p_i x_i)_i - (\subscr{P}{r}/n)\ones_n\big) + \theta \ones_n}{\theta\in\real} \\
	&= \setdef{-L^\dagger(p_i x_i)_i+\theta\ones_n}{\theta\in\real}.
	\end{aligned}
	\end{equation*}
	Moreover, substituting a $y^\star$ gives $\sigma
	\in \spn\{\ones_n\}$, and it follows that the problem $\P 2$
	reduces precisely to $\P 1$. \QED

\textbf{Proof of Lemma~\ref{lem:y-domain}:}
	The proof is trivially seen by multiplying $\dot{y}$ in~\eqref{eq:binpad} from the left by $\ones_n$ and applying the null space of $L$. \QED
	
\textbf{Proof of Lemma~\ref{lem:y-soln}:}
	The first term is computed by setting $\nabla_y\tilde{f}(x,y^\star) = \zeros_n$ (resp. $\nabla_y\widetilde{E}(x,y^\star) = \zeros_n)$ and solving for $y^\star$. There is a hyperplane of possible solutions due to the rank deficiency of $L$, but we are looking for the unique solution in $\Y$. The second term therefore follows from~\eqref{eq:ydom}.
	The fact that this point is also the unique equilibrium in $\Y$ follows from the fact that $\dot{y} = -\alpha\nabla_y \widetilde{E}(x,y^\star)$. \QED

\textbf{Proof of Lemma~\ref{lem:finite-eq-dist}:}
	The proof follows closely to the proof of Lemma~\ref{lem:finite-eq} with the variation that $\tilde{v}$ in the expression for $\dot{x}$ is now a function of $y$. Given the result of Lemma~\ref{lem:y-soln}, we may directly substitute the unique $y^\star$~\eqref{eq:ystar} for any $x$. Because $y^\star$ is simply a linear expression in $x$, the same argument as in Lemma~\ref{lem:finite-eq} that $\widetilde{\X}$ is finite follows. \QED
	
\textbf{Proof of Theorem~\ref{thm:dist-cvg}:}
	The first part of the proof to establish convergence to a critical point follows
	from a similar argument to the proof of
	Theorem~\ref{thm:cent-cvg}. Differentiating $\widetilde{E}$ with
	respect to time gives:
	\begin{equation}\label{eq:dEtildedt}
	\begin{aligned}
	\dfrac{d\widetilde{E}}{dt} &= \begin{bmatrix}
	\dot{x} \\ \dot{y}
	\end{bmatrix}^\top \begin{bmatrix}
	\nabla_x \widetilde{E}(x,y) \\ \nabla_y \widetilde{E}(x,y)
	\end{bmatrix} \\
	&= \begin{bmatrix}
	\dot{x} \\ \dot{y}
	\end{bmatrix}^\top 
	\begin{bmatrix}
	-\widetilde{W}x -\tilde{v} + g^{-1}(x)/\tau \\ -\alpha^{-1}\dot{y}
	\end{bmatrix} \\
	&= -\dot{x}^\top\diag{T/(x_i-x_i^2)_i}
	\vert \tilde{H}(x)\vert_m\dot{x} - \alpha^{-1}\dot{y}^\top \dot{y} < 0, \\
	& \qquad \qquad \qquad \dot{x} \neq 0 \ \text{or} \ \dot{y} \neq
	0, \quad (x,y)\in (0,1)^n \times \mathcal{Y}.
	\end{aligned}
	\end{equation}
	Thus, $\widetilde{E}$ monotonically decreases along the
	trajectories of $\binpad$. Given~\eqref{eq:dEtildedt}, we call again on the forward invariance property of the open hypercube for the distributed case via Lemma~\ref{lem:fwd-inv-d}, which verifies that $(x,y)\in (0,1)^n \times \mathcal{Y}$ at all times.

	Due to the deficiency induced by $L$, $\widetilde{E}$ is not radially unbounded in $y$ over all of $\real^n$, so we must be careful before applying the LaSalle Invariance Principle. Instead, define $\widetilde{E}$ only on $[0,1]^n\times \Y$ in consideration of Lemma~\ref{lem:y-domain}. Radial unboundedness in $\widetilde{E}$ is then obtained given any $y(0)$, and it follows that
	the trajectories converge to
	largest invariant set contained in $d\widetilde{E}/dt
	= 0$ per the LaSalle Invariance Principle~\cite{HKK:02}. This is the finite set of critical points of $\widetilde{E}$ per Lemma~\ref{lem:finite-eq-dist}, and so it additionally follows that we converge to a single critical point $(x^\star,y^\star)$.  
	
	Because $\widetilde{E}$ is convex in $y$, it follows that
	for any fixed $x$ there exist only local minima of
	$\widetilde{E}$ with respect to $y$. In consideration of this,
	we need only apply the Stable Manifold Theorem~\cite{JG-PH:83} to $x$. The
	argument for this develops similarly to the proof of
	Theorem~\ref{thm:cent-cvg}, and we conclude that the
	trajectories of $\binpad$ converge to a local minimizer
	$(x^\star,y^\star)$ of $\widetilde{E}$ with probability
	one. 
	
	The final part of the Theorem statement
	that $\widetilde{E}(x^\star,y) \geq \widetilde{E}(x^\star,y^\star),
	\forall y\in\real^n$ can also be seen from the convexity of $\widetilde{E}$
	in $y$ and applying the first-order condition of
	convexity:
	\begin{equation*}
	\widetilde{E}(x^\star,y) \geq \widetilde{E}(x^\star,y^\star) + (y-y^\star)^\top \nabla_y \widetilde{E}(x^\star,y^\star)
	\end{equation*}
	along with $\nabla_y\widetilde{E}(x^\star,y^\star) =
	\zeros_n$. \QED
	
\textbf{Proof of Lemma~\ref{lem:fwd-inv-d}:}
	The forward invariance of $\Y$ is already established per its definition and Lemma~\ref{lem:y-domain}, but we must establish that the trajectories $y(t)$ remain bounded in order to apply the argument in Lemma~\ref{lem:fwd-inv-c} to the proof of Theorem~\ref{thm:dist-cvg}. Compute the Hessian of $\widetilde{E}$ with respect to $y$ as:
	\begin{equation*}
	\nabla_{yy}\widetilde{E} = \gamma L^2 \succeq 0.
	\end{equation*}
	Due to the connectedness of $L$, the eigenspace associated with the $n-1$ strictly positive eigenvalues of $\gamma L^2$ is parallel to $\Y$. Therefore, $\widetilde{E}$ is strictly convex in $y$ on this subspace, and it follows that $\widetilde{E}$ is bounded from below on $\Y$. Due to $d\widetilde{E}/dt \leq 0$~\eqref{eq:dEtildedt} and the continuity of $\widetilde{E}$ in $y$, it follows that $y(t)$ is bounded for all $t$. Given this, the argument from Lemma~\ref{lem:fwd-inv-c} applies to the trajectories $x(t)$, and the set $(0,1)^n \times \Y$ is forward invariant under $\binpad$~\eqref{eq:binpad}. \QED

\bibliographystyle{abbrv}
\bibliography{alias,SMD-add,JC,SM}

\end{document}